\documentclass{amsart}

\usepackage{amsmath, amsfonts, amsthm, amssymb, dsfont}
\usepackage{exscale}
\usepackage{cite}
\usepackage{epsfig}
\usepackage{amscd}
\usepackage{graphics}
\usepackage{color}

\renewcommand{\geq}{\geqslant}

\newtheorem*{main-theorem}{Main Theorem}
\newtheorem*{theorem*}{Theorem}

\theoremstyle{definition}

\newtheorem*{remark*}{Remark}

\numberwithin{equation}{section}

\def\phi{\varphi}

\def\phi{\varphi}

\def\be{\begin{eqnarray*}}
\def\ee{\end{eqnarray*}}
\def\ben{\begin{eqnarray}}
\def\een{\end{eqnarray}}

\def\L2R{L_{\text{Rest}}^2}

\def\11{\mathds{1}}

\def\L2c{L^2_{\text{comp}}}

\usepackage{amsmath, amssymb, amsthm, enumitem}
\setlist{itemsep=.5em}

\newcommand{\mc}{\mathcal}

\newcommand{\begindescr}{\begin{description}[leftmargin=1.4em]}

\newcommand{\R}{\mathbb{R}}

\newcommand{\supinf}{^{\infty}}

\newcommand{\subk}{_{k=1}}

\newcommand{\de}{\,d}
\newcommand{\by}{\times}

\newcommand{\abs}[1]{\lvert#1\rvert}
\newcommand{\Abs}[1]{\left\lvert#1\right\rvert}
\newcommand{\norm}[2][]{\left\lVert#2\right\rVert_{#1}}
\newcommand{\paren}[1]{\left(#1\right)}
\newcommand{\brac}[1]{\left[#1\right]}

\newcommand{\defrac}[3][]{\frac{d^{#1}#2}{d{#3}^{#1}}}

\newcommand{\oline}[1]{\overline{#1}}

\DeclareMathOperator{\re}{Re}

\newtheorem{thm}{Theorem}
\theoremstyle{remark}
\newtheorem{rmk}{Remark}
\theoremstyle{lemma}
\newtheorem{lem}{Lemma}

\DeclareMathOperator{\area}{Area}
\title[Triangle observability]{Asymptotic Boundary Observability for the Wave Equation on One Side
  of a Planar Triangle}
\author{Hans Christianson}
\author{Evan Stafford}
\begin{document}
\begin{abstract}

We consider the wave equation $(\partial_t^2-\Delta)u=0$ on a planar triangular domain $\Omega\subset\R^2$ 
with Dirichlet boundary conditions. We use a commutator and
integration by parts argument similar to that in \cite{Chr2DTriangles}
by the first author 
to obtain an observability asymptotic 
for any one side of the triangle.  Our result is particular to triangular domains and does not hold for polygons in general.

  \end{abstract}

\maketitle

\section{Introduction}

In this paper, we study boundary observability for solutions to the
wave equation on flat triangles.  The main result is a large time asymptotic
observability identity from any one side.  In the closely related work 
\cite{Chr2DTriangles} by the first author,
it is shown that the Neumann data for Dirichlet eigenfunctions on a
triangle is equidistributed on each side.  Moreover, it is an exact
equality, so it is reasonable to expect a similar statement for
solutions to the wave equation.  We formulate this statement in terms
of an observability identity which roughly asks ``can we detect a wave
by making a measurement at the boundary?''
%% , which examines the Dirichlet eigenfunction problem on planar triangles and obtains the result that the $L^2$ norm of the Neumann data along any given side is equal to the length of the side divided by the area of the triangle. Here, we consider solutions to the wave equation on such triangular domains and proceed with the intention of obtaining an observability inequality. A noteworthy aspect of the result of \cite{Chr2DTriangles} is that it is an exact equality and not an estimate, and in this work we similarly are able to obtain a stronger result in the form of an asymptotic. 

The traditional methods for understanding control and observability
use geometric optics and microlocal analysis, however our approach
does not.  Instead, our result  only uses a robust commutator and
integration by parts argument. %% This requires some justification as triangular domains do not have smooth boundaries, but our use is nonetheless valid and is addressed in Appendix \ref{app_IBP}.
We begin the proof by computing the commutator of the wave operator
and an appropriate radial vector field. By using this computation in
tandem with evaluating the commutator explicitly, we obtain an
equation from which we derive our result, after applying a few
rudimentary 
estimates.  The paper is entirely self-contained.

 Let $\Omega$ be a triangular region in $\R^2$, and suppose the sides
  are labeled $A$, $B$, and $C$.  Let $\ell_A$ (respectively $\ell_B$, $\ell_C$)
  denote the altitude corresponding to $A$ (respectively $B, C$),
  meaning the perpendicular distance
  from the line containing $A, B, $ or $C$ to the non-adjacent corner
  (see Figures \ref{fig_acute} and \ref{fig_obtuse} for a picture).  Let $L$ denote the length of
  the longest side.  
Our main theorem is the following: 

\begin{thm}\label{thm_main}

  %% Let $\Omega$ be a triangular region in $\R^2$, and suppose the sides
  %% are labeled $F_1$, $F_2$, and $F_3$.  For $j = 1,2,3$ let $\ell_j$
  %% denote the altitude 

  %% Let $L$ be the length of the longest side.
 Consider the following initial/boundary value problem for the wave equation: 
	\begin{equation}\label{eq_main}
		\begin{cases}
			(\partial_t^2-\Delta)u=0 \\
			u|_{\partial\Omega} = 0 \\
			u(0,x,y) = u_0(x,y),\ u_t(0,x,y) = u_1(x,y)
		\end{cases}
	\end{equation}
	where $u_0, u_1\in H^1_0(\Omega) \cap H^s(\Omega)$ for all
        $s\ge0$. Then  for all $T>0$, the Neumann data on side $A$ satisfies 
	\begin{equation}\label{eq_result}
		\int_0^T\int_{A} \abs{\partial_\nu u}^2 \,dS \,dt =
                \frac{T}{\ell_A}E(0) \paren{1 + \mc{O} \paren{\frac{L}{T}}} 
	\end{equation}
	where $\partial_\nu$ is the normal derivative on $\partial\Omega$, $\de S$ is the arclength measure, 
	and $E(0)$ is the (conserved) initial energy, defined by 
	\begin{equation}\label{eq_energy}
		E(t) = \int_\Omega \abs{\partial_t u}^2 + \abs{\nabla u}^2 \de V.
	\end{equation}
        The analogous asymptotic on sides $B$ and $C$ also holds.
\end{thm}

%% \begin{rmk}
%% 	The $\mc{O}(1/T)$ term is calculated explicitly in the proof, namely in \eqref{eq_bound1} and \eqref{eq_bound2}. It is not optimal however.
%% \end{rmk}

\begin{rmk}
   This is not a
  low-regularity result; the assumptions on the regularity of $u_0$ and $u_1$ are simply to
  allow us to integrate by parts as necessary.
  \end{rmk}

\begin{rmk}
  Note that if $\partial_\nu u = 0$ on any one of the sides,
  we have that $E(0) = 0$, which implies $u=0$. In particular, given
  two solutions $u$ and $v$ to \eqref{eq_main} with identical Neumann
  data on at least one side, we have in fact that $u-v=0$ everywhere. Thus,
  Theorem \ref{thm_main} can be seen to imply a  uniqueness condition
  from Neumann data on just one side of the triangle.
\end{rmk}

\begin{rmk}
	It is important that the domain $\Omega$ is a triangle, as this theorem fails in general if 
	$\Omega$ is merely assumed to be polygonal. This is
        demonstrated in Section  \ref{ex_square}, where we show that
        the result fails for square domains.

        A similar result is believed to hold in higher dimensions and
        we will return to this question later.
\end{rmk}

\begin{rmk}
	The appearance of the factor $1/\ell_A$ in \eqref{eq_result}
        can be seen as a consequence of finite propagation speed,
since it takes time $\sim \ell_A$ for a wave to travel from the
opposite  corner
to side $A$.

        %% This is more easily seen for the one-dimensional analog of the problem, where $\Omega$ is an interval of length $\ell$ (this case is detailed in \S\ref{1D_case}). If $\Omega$ is the unit interval, a wave must travel at most a distance of 2 to meet the right endpoint. If $\Omega$ is scaled by a factor of $\ell$, then the bound on the distance a wave must travel to meet the right endpoint becomes $2\ell$. We can then expect the fraction of the energy that meets the right endpoint by a certain time to scale by a factor of $1/\ell$. 
\end{rmk}

\begin{rmk}[History] 
	For control and observability for solutions to the wave
        equation, a landmark result is one of Rauch-Taylor \cite{RT},
        who introduced the use of rays from geometric optics into the
        study of interior control problems.
        For boundary contol and observability as studied in this
        paper, an intuitive condition for
        the region of control to satisfy in order to establish
        observability is that every ray hits the region   in some
        finite time. This condition was shown to be necessary by a
        result of Ralston \cite{Ral69}, and \cite{RT} showed that for
        manifolds without boundary and interior control, this
        condition is (roughly) sufficient as well. Following this is an important result of Bardos-Lebeau-Rauch \cite{BLR}, which showed sufficiency for the case of observation and control problems on bounded domains for which the region of control lies on the boundary. 
	
	There are also results for cases where this geometric control
        condition is not fully satisfied. Lebeau \cite{Leb96}
        considered a particular example of a damped wave problem on a
        manifold where the damping region does not meet every
        geometric ray, but has enough control that the energy still
        decays exponentially with some derivative loss.  The papers
        \cite{Chr07} along with \cite{Chr10} by the first author built
        on this example and obtained a sub-exponential estimate that
        holds for a wider class of operators and manifolds. The paper \cite{CB15} demonstrated that this estimate is in general sharp, but can be improved with the addition of a stronger damping term.

\end{rmk}
\section{Simplified Cases}

\subsection{Case of a Single Dirichlet Eigenfunction}\label{ex_efun} %===== Single Eigenfunction Example ====

We first consider a case where we can make an explicit calculation. Let $\phi(x,y)$ be a Dirichlet eigenfunction on a triangle $\Omega$ satisfying 
\begin{equation}\label{eq_laplace}
	\begin{cases}
		(-\Delta - \lambda^2)\phi = 0 \\
		\phi|_{\partial\Omega} = 0,
	\end{cases}
\end{equation}
normalized so that $\norm[L^2(\Omega)]{\phi} = 1$. Note that by the main result of \cite{Chr2DTriangles}, we have that 
\begin{equation}\label{eq_efun_data}
	\int_A \abs{\partial_\nu \phi}^2 \de S = \frac{\lambda^2 a}{\area(\Omega)} = \frac{\lambda^2 a}{a\ell_A/2} 
	= \frac{2 \lambda^2}{\ell_A}
\end{equation}
where $a$ is the length of side $A$. Now define 
\begin{equation}
	u(t,x,y) = \sin(t\lambda)\phi(x,y),
\end{equation}
so $u$ satisfies 
\begin{equation}
	\begin{cases}
		(\partial_t^2 - \Delta)u=0 \\
		u|_{\partial\Omega} = 0 \\
		u(0,x,y) = 0,\ \partial_tu(0,x,y) = \lambda \phi(x,y).
	\end{cases}
\end{equation}
Note that for solutions to the wave equation \eqref{eq_main}, the energy \eqref{eq_energy} is conserved: 
\begin{align}
\begin{split}
	\defrac{}{t} E(t) &= \int_\Omega \partial_t \oline{u} \partial_t^2 u + \partial_t u \partial_t^2 \oline{u} 
			+ (\partial_t \nabla \oline{u})\cdot(\nabla u) + (\partial_t \nabla u)\cdot(\nabla \oline{u}) \de V \\ 
		&= \int_\Omega \partial_t \oline{u} \partial_t^2 u + \partial_t u \partial_t^2 \oline{u} 
			- \partial_t \oline{u}\Delta u - \partial_t u \Delta \oline{u} \de V \\
			&\qquad + \int_{\partial\Omega} \partial_t \oline{u}\partial_\nu u + \partial_t u \partial_\nu\oline{u}\de S \\
		&= \int_\Omega \partial_t\oline{u}(\partial_t^2 u - \Delta u) + \partial_t u (\partial_t^2\oline{u} - \Delta\oline{u}) \de V \\
		&= 0
\end{split}
\end{align}
by integration by parts and using the conditions in \eqref{eq_main}. In this case, since 
$\norm[L^2(\Omega)]{\phi}=1$, the energy is 
\begin{align}
\begin{split}
	E(0) &= \int_\Omega \abs{\partial_t u(0,x,y)}^2 + \abs{\nabla u(0,x,y)}^2 \de V \\
		&= \int_\Omega \lambda^2\abs{\phi}^2 \de V \\
		&= \lambda^2.
\end{split}
\end{align}
We can then compute explicitly 
\begin{align}
\begin{split}
	\int_0^T \int_A \abs{\partial_\nu u}^2 \de S \de t 
		&= \int_0^T \int_A \sin^2(t\lambda) \abs{\partial_\nu\phi}^2 \de S \de t \\
		&= \int_0^T \frac12 (1-\cos(2t\lambda)) \int_A \abs{\partial_\nu \phi}^2 \de S \de t \\
		&= \frac{\lambda^2}{\ell_A} \paren{T - \frac{1}{2\lambda} \sin(2T\lambda)} \\
		&= \frac{T}{\ell_A} E(0)\paren{1  - \frac{1}{2T\lambda} \sin(2T\lambda)}
\end{split}
\end{align}
and we see that the result is realized in this case. 

As Dirichlet eigenfunctions $\{\phi_k\}\subk\supinf$ form an orthonormal basis of $L^2(\Omega)$, we may 
ask whether our main result could be established by a similar calculation after noting that a solution $u$ to 
\eqref{eq_main} may be expressed as an expansion of the form 
\begin{equation}
	u(t,x,y) = \sum\subk\supinf \paren{a_k e^{it\lambda_k} + b_k e^{-it\lambda_k}} \phi_k(x,y).
\end{equation}
Unfortunately, at this time this is not possible. In the computation of $\int_A \abs{\partial_\nu u}^2 \de S$ we would in general obtain cross terms, and we do not have that the $\partial_\nu \phi_k|_{\partial\Omega}$ are orthogonal in $L^2(\partial\Omega)$. Thus, another approach is needed, which we will first demonstrate for the simpler one-dimensional case in \S\ref{1D_case} and later adapt to the two-dimensional case for the main proof in Section \ref{ch_proof}.

\subsection{One-Dimensional Problem}\label{1D_case} %===================== 1D Case =======

Here we present the one-dimensional analog of the problem stated in Theorem \ref{thm_main} in order to establish our general approach. 

\begin{thm}\label{thm_1D}
Let $\Omega\subset\R$ be the interval $[0,\ell]$, and consider the following initial/boundary value problem for the wave equation: 
\begin{equation}\label{1D_wave_eq}
	\begin{cases}
		(\partial_t^2 - \partial_x^2)u = 0 \\
		u(t,0) = u(t,\ell) = 0 \\
		u(0,x) = u_0(x),\ u_t(0,x) = u_1(x)
	\end{cases}
\end{equation}
where $u_0, u_1 \in H^1_0(\Omega) \cap H^s(\Omega)$ for all $s\ge0$. Then for all $T>0$, the Neumann data at the endpoint $x=\ell$ satisfies 
\begin{equation}
	\int_0^T \abs{\partial_x u}^2 \Big\vert_{x=\ell} \de t = \frac{T}{\ell} E(0) \paren{1+\mc{O}\paren{\frac{\ell}{T}}}
\end{equation}
where $E(0)$ is the (conserved) initial energy, defined by 
\begin{equation}
	E(t) = \int_0^\ell \abs{\partial_t u}^2 + \abs{\partial_x u}^2 \de x.
\end{equation}
\end{thm}

To see that the energy is conserved, we compute 
\begin{align}
\begin{split}
	\defrac{}{t} E(t) &= \defrac{}{t} \int_0^\ell \abs{\partial_t u}^2 + \abs{\partial_x u}^2 \de x \\
		&= \int_0^\ell \partial_t^2 u \partial_t \oline{u} + \partial_t u \partial_t^2 \oline{u} 
			+ \partial_x\partial_t u \partial_x\oline{u} + \partial_x u \partial_x\partial_t \oline{u} \de x \\
		&= \int_0^\ell \partial_t^2 u \partial_t \oline{u} + \partial_t u \partial_t^2 \oline{u} 
			- \partial_t u \partial_x^2\oline{u} - \partial_x^2 u \partial_t \oline{u} \de x 
			+ \brac{\partial_t u \partial_x \oline{u} + \partial_x u \partial_t \oline{u}} \Big\vert_0^\ell \\
		&= \int_0^\ell (\partial_t^2 u - \partial_x^2 u) \partial_t\oline{u} 
			+ (\partial_t^2\oline{u} - \partial_x^2\oline{u}) \partial_t u \de x \\
		&= 0
\end{split}
\end{align}
by integration by parts and using the conditions in \eqref{1D_wave_eq}. 

Consider the vector field 
\begin{equation}
	X = x\partial_x
\end{equation}
on $\Omega$. We have the commutator
\begin{align}
\begin{split}
	[\partial_t^2 - \partial_x^2, X] 
		&= (\partial_t^2 - \partial_x^2) x \partial_x - x\partial_x(\partial_t^2-\partial_x^2) \\
		&= -\partial_x^2 x \partial_x + x \partial_x^3 \\
		&= -\partial_x (\partial_x + x \partial_x^2) + x \partial_x^3 \\
		&= -\partial_x^2 - \partial_x^2 - x\partial_x^3 + x \partial_x^3 \\
		&= -2\partial_x^2.
\end{split}
\end{align}
Using this, along with integration by parts and the fact that $u$ satisfies the wave equation, we obtain 
\begin{align}\label{comm_eq_start}
\begin{split}
	\int_0^T \int_0^\ell [\partial_t^2 - \partial_x^2, X] u\oline{u} \de x \de t 
		&= \int_0^T \int_0^\ell -2\partial_x^2u \oline{u} \de x \de t \\
		&= \int_0^T \int_0^\ell (-\partial_t^2-\partial_x^2)u \oline{u} \de x \de t \\
		&= \int_0^T \int_0^\ell \abs{\partial_tu}^2 + \abs{\partial_xu}^2 \de x \de t \\
		&\qquad - \int_0^\ell \partial_t u\oline{u} \de x \Big\vert_0^T - \int_0^T \partial_x u\oline{u} \Big\vert_0^\ell \de t. 
\end{split}
\end{align}
We recognize the energy in the first term, and we note that the third term vanishes by the boundary conditions in \eqref{1D_wave_eq}, so 
\begin{align}\label{comm_eq_rhs}
\begin{split}
	\int_0^T \int_0^\ell [\partial_t^2 - \partial_x^2, X] u\oline{u} \de x \de t 
		&= \int_0^T E(0) \de t - \int_0^\ell \partial_t u\oline{u} \de x \Big\vert_0^T \\
		&= TE(0) - \int_0^\ell \partial_t u\oline{u} \de x \Big\vert_0^T. 
\end{split}
\end{align}

From the starting point of \eqref{comm_eq_start}, we now evaluate the commutator explicitly, which yields 
\begin{align}
\begin{split}
	\int_0^T \int_0^\ell [\partial_t^2 - \partial_x^2, X] u\oline{u} \de x \de t 
		&= \int_0^T \int_0^\ell (\partial_t^2 - \partial_x^2) x\partial_x u\oline{u} 
			- x\partial_x(\partial_t^2 - \partial_x^2) u\oline{u} \de x \de t \\ 
		&= \int_0^T \int_0^\ell \partial_t^2 x\partial_x u\oline{u} 
			- \partial_x^2 x\partial_x u\oline{u} \de x \de t \\
		&= \int_0^T \int_0^\ell - (\partial_t (x \partial_x u))(\partial_t u) 
			+ (\partial_x (x \partial_x u))(\partial_x\oline{u}) \de x \de t \\
		&\qquad + \int_0^\ell (\partial_t(x\partial_x u))\oline{u} \de x \Big\vert_0^T
			- \int_0^T (\partial_x(x\partial_xu))\oline{u} \Big\vert_0^\ell \de t
\end{split}
\end{align}
again by integration by parts and the fact that $u$ satisfies the wave equation. Note that the third term vanishes by the boundary conditions. We may integrate by parts in the first term again and use the fact that $\oline{u}$ also satisfies the wave equation to obtain 
\begin{align}
\begin{split}
	\int_0^T \int_0^\ell [\partial_t^2 - \partial_x^2, X] u\oline{u} \de x \de t 
		&= \int_0^T \int_0^\ell (x\partial_x u)(\partial_t^2\oline{u}) 
			- (x\partial_xu)(\partial_x^2\oline{u}) \de x \de t \\
		&\qquad + \int_0^\ell (\partial_t(x\partial_x u))\oline{u} 
			- (x\partial_x u)(\partial_t \oline{u}) \de x \Big\vert_0^T \\
		&\qquad + \int_0^T (x\partial_x u)(\partial_x\oline{u})\Big\vert_0^\ell \de t \\
		&= \int_0^\ell (\partial_t(x\partial_x u))\oline{u} 
			- (x\partial_x u)(\partial_t \oline{u}) \de x\Big\vert_0^T \\
		&\qquad + \int_0^T \ell\abs{\partial_x u}^2 \Big\vert_{x=\ell} \de t.
\end{split}
\end{align}
Using the fact that $x\partial_xu = \partial_x(xu)-u$, we can then integrate by parts yet again to arrive at 
\begin{align}\label{comm_explicit}
\begin{split}
	\int_0^T \int_0^\ell [\partial_t^2 - \partial_x^2, X] u\oline{u} \de x \de t 
		&= \int_0^\ell (\partial_t(\partial_x(xu))\oline{u} - \partial_t u \oline{u} 
			- (x\partial_x u)(\partial_t \oline{u}) \de x \Big\vert_0^T \\
		&\qquad + \int_0^T \ell\abs{\partial_x u}^2 \Big\vert_{x=\ell} \de t \\
		&= \int_0^\ell -(\partial_tu)(x\partial_x\oline{u}) - \partial_t u \oline{u} 
			- (x\partial_x u)(\partial_t \oline{u}) \de x \Big\vert_0^T \\
		&\qquad + \int_0^T \ell\abs{\partial_x u}^2 \Big\vert_{x=\ell} \de t 
			+ \partial_t(xu)\oline{u} \Big\vert_{x=0}^\ell \Big\vert_{t=0}^T \\
		&= \int_0^\ell -(\partial_tu)(x\partial_x\oline{u}) - \partial_t u \oline{u} 
			- (x\partial_x u)(\partial_t \oline{u}) \de x \Big\vert_0^T \\
		&\qquad + \int_0^T \ell\abs{\partial_x u}^2 \Big\vert_{x=\ell} \de t.
\end{split}
\end{align}
Combining this with \eqref{comm_eq_rhs} and reorganizing terms then yields
\begin{equation}\label{1D_comm_eq}
	\int_0^T \ell\abs{\partial_x u}^2 \Big\vert_{x=\ell} \de t = TE(0) 
		+ \int_0^\ell (\partial_tu)(x\partial_x\oline{u}) + (x\partial_x u)(\partial_t \oline{u}) \de x \Big\vert_0^T.
\end{equation}
Note that in this case, the $\int_0^\ell \partial_tu\oline{u} \de x \big\vert_0^T$ terms from \eqref{comm_eq_rhs} and \eqref{comm_explicit} cancel. This does not occur in the two-dimensional case, where these extra terms will instead require us to make additional estimates to control them.

We now obtain bounds for the integral on the right hand side of \eqref{1D_comm_eq}. By Cauchy's inequality with parameter $\alpha>0$  we have 
\begin{align}
\begin{split}
	\int_0^\ell \abs{(\partial_tu)(x\partial_x\oline{u})} \de x 
		&\le \int_0^\ell \frac\alpha2 \abs{\partial_t u}^2 + \frac{1}{2\alpha} \abs{x\partial_x u}^2 \de x \\
		&\le \int_0^\ell \frac\alpha2 \abs{\partial_t u}^2 + \frac{\ell^2}{2\alpha} \abs{\partial_x u}^2 \de x.
\end{split}
\end{align}
Taking $\alpha = \ell$, this becomes 
\begin{equation}
	\int_0^\ell \abs{(\partial_tu)(x\partial_x\oline{u})} \de x \le \frac{\ell}{2} E(0).
\end{equation}
Equation \eqref{1D_comm_eq} then implies simultaneously 
\begin{equation}
	\int_0^T \ell\abs{\partial_x u}^2 \Big\vert_{x=\ell} \de t \le TE(0) + 4\frac{\ell}{2} E(0)
\end{equation}
and 
\begin{equation}
	\int_0^T \ell\abs{\partial_x u}^2 \Big\vert_{x=\ell} \de t \ge TE(0) - 4\frac{\ell}{2} E(0),
\end{equation}
which yield the asymptotic 
\begin{equation}
	\int_0^T \abs{\partial_x u}^2 \Big\vert_{x=\ell} \de t = \frac{T}{\ell} E(0) \paren{1+\mc{O}\paren{\frac{\ell}{T}}}.
\end{equation}
This is the one-dimensional analog of the main result. We now move on
to the proof of the two dimensional case.

\section{Proof of the Main Theorem}\label{ch_proof}

\subsection{Proof for Acute Triangles} %====================== Acute Proof ======================

The proof of the main theorem is split into two cases, one for acute
triangles and one for obtuse (or right) triangles. The steps in both
cases are the same in spirit, and in fact the proof from the acute
case works for all but one side in the obtuse case.  Without loss in generality, we prove
the main theorem for side $A$ only, and use $\ell = \ell_A$ to ease
notation.

    \begin{figure}
\hfill
\centerline{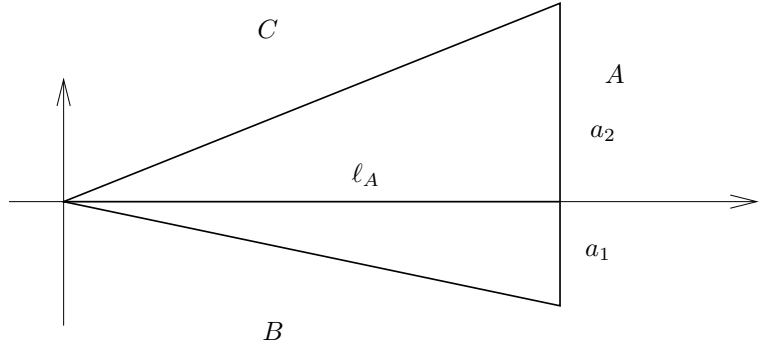}
\caption{\label{fig_acute}  Setup for acute  triangles.}
\hfill
\end{figure}

%% \begin{figure}[h]
%% 	\begin{center}
%% 		\begin{tikzpicture}[auto, decoration={markings,mark=at position 0.55 with {\arrow[line width=.5mm]{>}}}, scale=1, thick]
%% 		\draw[<->] (-2,0) to (9,0);
%% 		\draw[<->] (0,4) to (0,-3);
		
%% 		\draw[thick] (0,0) to node {$C$} (7,3);
%% 		\draw[thick] (0,0) to node[swap] {$B$} (7,-2); 
%% 		\draw[thick] (7,3) to node {$A$} (7,-2);
		
%% 		\node at (5, 1) {$\Omega$};
%% 		\node at (4, -0.3) {$\ell$};
%% 		\node at (6.6,1.3) {$a_2$};
%% 		\node at (6.6,-.8) {$a_1$};
		
%% 		\end{tikzpicture}
%% 	\end{center}
%% 	\caption{Acute triangle case}\label{fig_acute}
%% \end{figure}

We use a setup similar to that in \cite{Chr2DTriangles}. Suppose $\Omega$ is an acute triangle, as in Figure \ref{fig_acute}. Let $\Omega$ be oriented such that the vertex opposing $A$ is located at the origin, and $A$ is perpendicular to the $x$-axis. Note that this causes the altitude of length $\ell$ to coincide with the $x$-axis. Label the remaining sides $B$ and $C$ in a clockwise fashion. Let $a_1$ be the length of the part of $A$ below the $x$-axis, and similarly let $a_2$ be the length of the part above. Consider the problem \eqref{eq_main} and define the vector field 
\begin{equation}\label{eq_vector_field_def}
	X = x\partial_x + y\partial_y
\end{equation}
on $\Omega$. We have the commutator 
\begin{align}\label{commutator}
\begin{split}
	[\partial_t^2-\Delta,X] &= (\partial_t^2-\Delta)X - X(\partial_t^2-\Delta) \\
		&= -\Delta x\partial_x - \Delta y\partial_y + X\Delta \\
		&= -2\partial_x^2 -x\partial_x\Delta - 2\partial_y^2 - y\partial_y\Delta + X\Delta \\
		&= -2\Delta.
\end{split}
\end{align}
From this and the fact that $u$ satisfies the wave equation, by integration by parts we have  
\begin{align}\label{eq_comm_lhs}
\begin{split}
	\int_0^T \int_\Omega [\partial_t^2-\Delta, X]u\overline{u} \,dV\,dt 
		&= \int_0^T \int_\Omega (-2\Delta u)\overline{u} \,dV\,dt \\
		&= \int_0^T \int_\Omega (-\Delta u - \partial_t^2u)\overline{u} \,dV\,dt \\
		&= \int_0^T \int_\Omega \lvert\nabla u\rvert^2 + \lvert\partial_t u\rvert^2 \,dV\,dt  
			+ \int_\Omega -\partial_t u \overline{u} \,dV \Big\vert_0^T \\
		&= TE(0) + \int_\Omega -\partial_t u \overline{u} \,dV \Big\vert_0^T, 
\end{split}
\end{align}
mirroring what was obtained in the one-dimensional case. Also, by evaluating the commutator explicitly and integrating by parts, we obtain 
\begin{align}
\begin{split}
	\int_0^T \int_\Omega [\partial_t^2-\Delta, X]u\overline{u} \,dV\,dt 
		&= \int_0^T \int_\Omega (\partial_t^2-\Delta) X u\overline{u} 
			- X(\partial_t^2-\Delta )u\overline{u} \,dV\,dt \\
		&= \int_0^T \int_\Omega \partial_t^2 X u\overline{u} -\Delta Xu\overline{u} \,dV\,dt \\
		&= \int_0^T \int_\Omega - (\partial_t Xu)(\partial_t\oline{u}) + (\nabla Xu)\cdot(\nabla\oline{u}) \,dV\,dt \\
		&\qquad + \int_\Omega (\partial_t Xu)\overline{u} \de V\Big\vert_0^T
			- \int_0^T\int_{\partial\Omega} (\partial_\nu Xu)\overline{u} \,dS \,dt. 
\end{split}
\end{align}
The third term vanishes since $\oline{u} |_{\partial\Omega} = 0$. Integrating by parts again in the first term then gives 
\begin{align}
\begin{split}
	\int_0^T \int_\Omega [\partial_t^2-\Delta, X]u\overline{u} \,dV\,dt 
		&= \int_0^T \int_\Omega (Xu)(\partial_t^2\overline{u}) - (Xu)(\Delta\overline{u}) \,dV\,dt \\
			&\qquad + \int_\Omega (\partial_t Xu)\overline{u} - (Xu)(\partial_t\overline{u}) \,dV \Big\vert_0^T \\
			&\qquad + \int_0^T\int_{\partial\Omega} (Xu)(\partial_\nu\overline{u}) \,dS \,dt.
\end{split}
\end{align}
Since $\oline{u}$ also satisfies the wave equation, the first term on the right disappears. Using the fact that 
$Xu = -2u + \partial_x(xu) + \partial_y(yu)$ and integrating by parts in the second term (the boundary term vanishes by the boundary conditions) then yields 
\begin{align}
\begin{split}
	\int_0^T \int_\Omega [\partial_t^2-\Delta, X]u\overline{u} \,dV\,dt 
		&= \int_\Omega -2\partial_t u \overline{u} -(\partial_t u)(X\overline{u}) - (Xu)(\partial_t\overline{u}) \,dV \Big\vert_0^T \\
		&\qquad + \int_0^T\int_{\partial\Omega} (Xu)(\partial_\nu\overline{u}) \,dS \,dt. 
\end{split}
\end{align}
Note that this is slightly different from what was obtained in the one-dimensional case due to the factor of 2 in the first term on the right hand side. Combining this with \eqref{eq_comm_lhs} and reorganizing terms gives 
\begin{equation}\label{eq_comm_equality}
	\int_0^T\int_{\partial\Omega} (Xu)(\partial_\nu\overline{u}) \,dS \,dt 
		= TE(0) + \int_\Omega \partial_t u \overline{u} + (\partial_t u)(X\overline{u}) + (Xu)(\partial_t\overline{u}) \,dV \Big\vert_0^T,
\end{equation}
and we see that there is a new term on the right which we must also account for.

First, we can obtain bounds in terms of the energy $E(0)$ for the
latter terms on the right hand side in a similar fashion to the
one-dimensional case. Note that $\abs{x}$ and $\abs{y}$ are bounded by
$L$ ($L = $ length of longest side as in the introduction) on $\Omega$ due to our choice of orientation for the triangle. An application of Cauchy's inequality with parameter $\alpha>0$ and the triangle inequality then implies 
\begin{align}
\begin{split}
	\int_\Omega \left| (\partial_t u) X \overline{u} \right| \,dV 
		&\le \int_\Omega \frac \alpha2 \abs{\partial_t u}^2 + \frac{1}{2\alpha} \abs{Xu}^2 \,dV \\
		&\le \int_\Omega \frac {\alpha}{2} \abs{\partial_t u}^2 + \frac{L^2}{2\alpha} \paren{\abs{\partial_x u} 
			+ \abs{\partial_y u}}^2 \,dV \\
		&= \int_\Omega \frac {\alpha}{2} \abs{\partial_t u}^2 + \frac{L^2}{2\alpha} \abs{\nabla u}^2 + \frac{L^2}{\alpha} \abs{\partial_x u}\abs{\partial_y u} \,dV
\end{split}
\end{align}
Applying Cauchy's inequality again then implies 
\begin{align}\label{eq_bound1}
\begin{split}
	\int_\Omega \abs{ (\partial_t u) (X \overline{u}) } \,dV 
		&\le \int_\Omega \frac{\alpha}{2} \abs{\partial_t u}^2 + \frac{L^2}{2\alpha} \abs{\nabla u}^2 + \frac{L^2}{2\alpha} \abs{\nabla u}^2 \,dV \\
		&\le \frac{L}{\sqrt{2}} E(0)
\end{split}
\end{align}
where $\alpha = L\sqrt{2}$ is chosen. 

For the new term, note that by Cauchy's inequality with parameter $\beta>0$, 
\begin{align}\label{eq_pre_poincare}
\begin{split}
	\Abs{ \int_\Omega \partial_t u \oline{u} \de V } 
		&\le \norm[L^2(\Omega)]{\partial_t u}\norm[L^2(\Omega)]{u} \\
		&\le \frac{\beta}{2} \norm[L^2(\Omega)]{\partial_t u}^2 + \frac{1}{2\beta} \norm[L^2(\Omega)]{u}^2.
\end{split}
\end{align}
We now establish a Poincar\'e-type inequality.  We present a full
proof here since it is very simple in the case of a triangle.
\begin{lem} Let $u \in H_0^1( \Omega ) \cap H^s ( \Omega)$ for every
  $s$.  Then the following holds:
	\begin{equation}\label{eq_poincare}
	\norm[L^2(\Omega)]{u} \le L\sqrt{e-1}\norm[L^2(\Omega)]{\partial_x u}.
\end{equation}
\end{lem}
We observe that the right hand side involves only the $x$ derivative
of $u$.  
This allows us to bound the right hand side of \eqref{eq_pre_poincare} using the energy after an appropriate choice of $\beta$. Such inequalities have been widely studied as described in \cite{KN15}, which notes results such as Poincar\'e \cite{Poincare90} and Steklov \cite{Steklov1897}. We note that if we happen to know the first Dirichlet eigenvalue $\lambda_1^2$ for \eqref{eq_laplace}, we readily have the inequality $\norm[L^2(\Omega)]{u} \le \lambda_1^{-2} \norm[L^2(\Omega)]{\nabla u}$ by Rayleigh's formula. Furthermore, for convex domains a result of Protter \cite{Protter81} demonstrates a lower bound on $\lambda_1^2$ in terms of the diameter of $\Omega$
and the radius of the largest disk contained in $\Omega$, so we could
then obtain \eqref{eq_poincare} with a constant that depends only on
$\Omega$, however it would involve both spatial derivatives. Here, we will derive the constant explicitly and see that it depends on $L$. 

\begin{proof}
Fix $t_0 \ge 0$ and define
\begin{equation}\label{eq_F}
	F(x) = \int_{-\frac{a_1}{\ell} x}^{\frac{a_2}{\ell} x} \abs{u(t_0,x,y)}^2 \de y,
\end{equation}
noting that 
\begin{equation}
	\int_0^\ell F(x) \de x = \norm[L^2(\Omega)]{u}^2. 
\end{equation}
Then by Cauchy's inequality with parameter $\gamma>0$ and the fact that $u(t_0,x,y)|_{\partial\Omega} = 0$, 
\begin{align}
\begin{split}
	F'(x) &= \frac{a_2}{\ell} \Abs{u\paren{t_0,x,\frac{a_2}{\ell} x}}^2 + \frac{a_1}{\ell} \Abs{u\paren{t_0,x,-\frac{a_1}{\ell} x}}^2 + 2 \re \int_{-\frac{a_1}{\ell} x}^{\frac{a_2}{\ell} x} u \partial_x \oline{u} \de y \\
		&\le \gamma \int_{-\frac{a_1}{\ell} x}^{\frac{a_2}{\ell} x} \abs{u}^2 \de y 
			+ \frac 1\gamma \int_{-\frac{a_1}{\ell} x}^{\frac{a_2}{\ell} x} \abs{\partial_x u}^2 \de y \\
		&= \gamma F(x) + \frac 1\gamma G(x)
\end{split}
\end{align}
where 
\begin{equation}\label{eq_G}
G(x) = \int_{-\frac{a_1}{\ell} x}^{\frac{a_2}{\ell} x} \abs{\partial_x u(t_0,x,y)}^2 \de y. 
\end{equation}
Then we have 
\begin{equation}
	F'(x) e^{-\gamma x} - \gamma F(x) e^{-\gamma x} \le \frac1\gamma e^{-\gamma x} G(x)
\end{equation}
and in particular,
\begin{equation}
	(F(x) e^{-\gamma x})' \le \frac{e^{-\gamma x}}{\gamma} G(x). 
\end{equation}
By an argument following a proof of Gronwall's inequality, we see that 
\begin{align}
\begin{split}
	F(x) &\le e^{\gamma x} \paren{F(0) + \int_0^x \frac{e^{-\gamma s}}{\gamma} G(s) \de s} \\
		&= \frac{e^{\gamma x}}{\gamma} \int_0^x e^{-\gamma s} G(s) \de s \\
		&\le \frac{e^{\gamma x}}{\gamma} \int_0^\ell \int_{-\frac{a_1}{\ell} s}^{\frac{a_2}{\ell} s} \abs{\partial_x u}^2 \de y \de s \\
		&= \frac{e^{\gamma x}}{\gamma} \norm[L^2(\Omega)]{\partial_x u}^2 
\end{split}
\end{align}
and we then obtain 
\begin{align}\label{eq_poincare_explicit}
\begin{split}
	\norm[L^2(\Omega)]{u}^2 &= \int_0^\ell F(x) \de x \\
		&\le \frac{1}{\gamma^2}\paren{e^{\gamma \ell} -1} \norm[L^2(\Omega)]{\partial_x u}^2.
\end{split}
\end{align}
Choosing $\gamma = \frac 1\ell$ then yields
\begin{equation}
	\norm[L^2(\Omega)]{u}^2 \le \ell^2\paren{e -1}
        \norm[L^2(\Omega)]{\partial_x u}^2 \le L^2\paren{e -1}
        \norm[L^2(\Omega)]{\partial_x u}^2
\end{equation}
since $L$ is the length of the longest side.  This 
 proves the lemma. This choice of $\gamma$ is probably not optimal, but it is sufficient for our purposes.
\end{proof}

Combining this with \eqref{eq_pre_poincare}, we then have 
\begin{align}
\begin{split}
	\Abs{ \int_\Omega \partial_t u \oline{u} \de V } &\le \frac{\beta}{2} \norm[L^2(\Omega)]{\partial_t u}^2 + \frac{1}{2\beta} \norm[L^2(\Omega)]{u}^2 \\
		&\le \frac{\beta}{2} \norm[L^2(\Omega)]{\partial_t
          u}^2 + \frac{L^2}{2\beta} \paren{e -1}
        \norm[L^2(\Omega)]{\partial_x u}^2 \\
      &\le \frac{\beta}{2} \norm[L^2(\Omega)]{\partial_t
          u}^2 + \frac{L^2}{2\beta} \paren{e -1}
        \norm[L^2(\Omega)]{\nabla u}^2  
\end{split}
\end{align}
which gives
\begin{equation}\label{eq_bound2}
	\Abs{ \int_\Omega \partial_t u \oline{u} \de V } \le \frac12 L \sqrt{e-1} E(0)
\end{equation}
after choosing $\beta = L\sqrt{e-1}$.

With \eqref{eq_bound1} and \eqref{eq_bound2} in hand, we then obtain from \eqref{eq_comm_equality} simultaneously 
\begin{equation}
	\Abs{\int_0^T\int_{\partial\Omega} (Xu)(\partial_\nu\overline{u}) \,dS \,dt} 
		\le TE(0) + 2\paren{2\frac{L}{\sqrt2}E(0) + \frac12 L \sqrt{e-1} E(0)}
\end{equation}
and
\begin{equation}
	\Abs{\int_0^T\int_{\partial\Omega} (Xu)(\partial_\nu\overline{u}) \,dS \,dt} 
		\ge TE(0) - 2\paren{2\frac{L}{\sqrt2}E(0) + \frac12 L \sqrt{e-1} E(0)},
\end{equation}
which we express together as the asymptotic 
\begin{equation}\label{ineq_with_rhs}
	\Abs{\int_0^T\int_{\partial\Omega} (Xu)(\partial_\nu\overline{u}) \,dS \,dt} 
		= TE(0)\paren{1 + \mc{O}\paren{\frac{L}{T}}}. 
\end{equation}

In the same fashion as in \cite{Chr2DTriangles}, we now get the Neumann data on side $A$ from the term on the left hand side. Since $u|_{\partial\Omega} = 0$, the tangential derivative of $u$ along each side vanishes. On $A$, this implies $\partial_y u = 0$. Since also $x=\ell$ and $\partial_\nu u = \partial_x u$ here, we have $Xu = \ell \partial_\nu u$, so 
\begin{align}\label{eq_A_data}
\begin{split}
	\int_0^T\int_A (Xu)(\partial_\nu\overline{u}) \,dS \,dt 
		&= \int_0^T\int_A \ell (\partial_\nu u)(\partial_\nu\overline{u}) \,dS \,dt \\
		&= \int_0^T\int_A \ell \abs{\partial_\nu u}^2 \,dS \,dt.
\end{split}
\end{align}

On side $B$, we have $y = -\frac{a_1}{\ell} x$, and the unit tangent vector is 
$\tau = \paren{\frac{\ell}{b}, -\frac{a_1}{b}}$. Since the tangential derivative vanishes, we have 
\begin{equation}\label{eq_B_relation}
	\tau \cdot \nabla u = \frac{\ell}{b} \partial_x u - \frac{a_1}{b} \partial_y u = 0.
\end{equation}
The unit normal vector on $B$ is $\nu = \paren{-\frac{a_1}{b}, -\frac{\ell}{b}}$. Using this and \eqref{eq_B_relation}, we obtain 
\begin{align}
\begin{split}
	(Xu)(\partial_\nu\oline{u}) &= \paren{x\partial_x u - \frac{a_1}{\ell}x\partial_y u} \paren{-\frac{a_1}{b}\partial_x\oline{u} - \frac{\ell}{b}\partial_y\oline{u}} \\
		&= x\paren{-\frac{a_1}{b}\partial_xu\partial_x\oline{u} + \frac{a_1^2}{b\ell} \partial_yu\partial_x\oline{u} 
			- \frac{\ell}{b}\partial_xu\partial_y\oline{u} + \frac{a_1}{b}\partial_yu\partial_y\oline{u}} \\
		&= -x\frac{a_1}{\ell}\partial_x\oline{u}\paren{\frac{\ell}{b} \partial_x u - \frac{a_1}{b} \partial_y u}
			- x\partial_y\oline{u}\paren{\frac{\ell}{b} \partial_x u - \frac{a_1}{b} \partial_y u} \\
		&= 0,
\end{split}
\end{align}
and consequently 
\begin{equation}\label{eq_B_data}
	\int_0^T\int_B (Xu)(\partial_\nu\overline{u}) \,dS \,dt = 0.
\end{equation}

In a similar fashion, on side $C$ we have $y = \frac{a_2}{c}x$, and the unit tangent vector is 
$\tau = \paren{\frac{\ell}{c}, \frac{a_2}{c}}$. Since the tangential derivative vanishes, we have 
\begin{equation}\label{eq_C_relation}
	\tau \cdot \nabla u = \frac{\ell}{c} \partial_x u + \frac{a_2}{c} \partial_y u = 0.
\end{equation}
The unit normal vector on $C$ is $\nu = \paren{-\frac{a_2}{c}, \frac{\ell}{c}}$. Using this and \eqref{eq_C_relation}, 
we obtain 
\begin{align}
\begin{split}
	(Xu)(\partial_\nu\oline{u}) &= \paren{x\partial_x u + \frac{a_2}{\ell}x\partial_y u} \paren{-\frac{a_2}{c}\partial_x\oline{u} + \frac{\ell}{c}\partial_y\oline{u}} \\
		&= x\paren{-\frac{a_2}{c}\partial_xu\partial_x\oline{u} - \frac{a_2^2}{c\ell} \partial_yu\partial_x\oline{u} 
			+ \frac{\ell}{c}\partial_xu\partial_y\oline{u} + \frac{a_2}{c}\partial_yu\partial_y\oline{u}} \\
		&= -x\frac{a_2}{\ell}\partial_x\oline{u}\paren{\frac{\ell}{c} \partial_x u + \frac{a_2}{c} \partial_y u}
			+ x\partial_y\oline{u}\paren{\frac{\ell}{c} \partial_x u + \frac{a_2}{c} \partial_y u} \\
		&= 0,
\end{split}
\end{align}
and consequently 
\begin{equation}\label{eq_C_data}
	\int_0^T\int_C (Xu)(\partial_\nu\overline{u}) \,dS \,dt = 0.
\end{equation}

Thus, from \eqref{eq_B_data} and \eqref{eq_C_data} we have 
\begin{equation}
	\int_0^T\int_{\partial\Omega} (Xu)(\partial_\nu\overline{u}) \,dS \,dt = \int_0^T\int_A \ell \abs{\partial_\nu u}^2 \,dS \,dt
\end{equation}
and \eqref{ineq_with_rhs} becomes 
\begin{equation}
	\int_0^T\int_A \ell \abs{\partial_\nu u}^2 \,dS \,dt = TE(0) \paren{1 + \mc{O}\paren{\frac{L}{T}}} ,
\end{equation}
implying 
\begin{equation}
	\int_0^T\int_A \abs{\partial_\nu u}^2 \,dS \,dt = \frac{T}{\ell}E(0) \paren{1 + \mc{O}\paren{\frac{L}{T}}} 
\end{equation}
which is the desired result.

\subsection{Proof for Obtuse (or Right) Triangles} %====================== Obtuse Proof ======================

    \begin{figure}
\hfill
\centerline{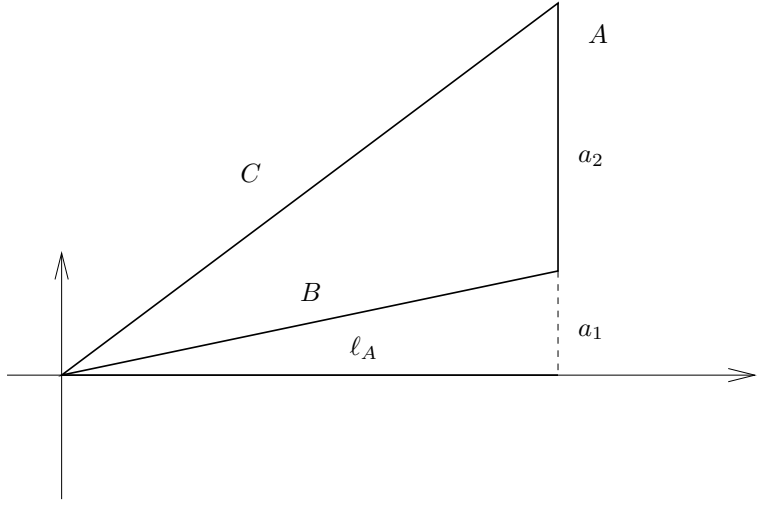}
\caption{\label{fig_obtuse}  Setup for obtuse (and right) triangles}
\hfill
\end{figure}

%% \begin{figure}[h]
%% \begin{center}
%% 	\begin{tikzpicture}[auto, decoration={markings,mark=at position 0.55 with {\arrow[line width=.5mm]{>}}}, scale=1, thick]
%% 	\draw[<->] (-2,0) to (9,0);
%% 	\draw[<->] (0,6) to (0,-1);
	
%% 	\draw[thick] (0,0) to node {$C$} (7,5);
%% 	\draw[thick] (0,0) to node[swap] {$B$} (7,1.5); 
%% 	\draw[thick] (7,5) to node {$A$} (7,1.5);
	
%% 	\draw[thick,dashed] (7,0) to (7,1.5);
	
%% 	\node at (5, 2) {$\Omega$};
%% 	\node at (4, -0.3) {$\ell$};
%% 	\node at (6.6,3) {$a_2$};
%% 	\node at (6.6,.75) {$a_1$};
	
%% 	\end{tikzpicture}
%% \end{center}
%% \caption{Obtuse (or right) triangle case}\label{fig_obtuse}
%% \end{figure}

The case for an obtuse triangle plays out very similarly. Let $\Omega$
now be an obtuse (or right) triangular region in $\R^2$, as in Figure
\ref{fig_obtuse}.
Before we begin, let us observe that the proof in the acute case
actually applies as long as $a_1 \geq 0$ in Figure \ref{fig_acute}.
That means that in Figure \ref{fig_obtuse} we only need to prove the
theorem for the side $A$ as the proof for sides $B$ and $C$ follow
from the acute case.

%% Since $A$ is chosen to correspond to the longest altitude, $A$ must be the shortest side, so the opposing angle will always be acute. Thus, the obtuse (or right) angle will not be located at the origin in our setup. 

Using the same vector field $X$ defined in \eqref{eq_vector_field_def}, we have the same commutator $[\partial_t^2 - \Delta, X] = -2\Delta$ as in \eqref{commutator}. The change in domain does not affect any of the steps leading up to equation \eqref{eq_comm_equality}, so it still holds in this case: 
\begin{equation}
	\int_0^T\int_{\partial\Omega} (Xu)(\partial_\nu\overline{u}) \,dS \,dt 
		= TE(0) + \int_\Omega \partial_t u \overline{u} + (\partial_t u)(X\overline{u}) + (Xu)(\partial_t\overline{u}) \,dV \Big\vert_0^T.
\end{equation}

We can obtain the same bounds as in \eqref{eq_bound1} for the last two terms, so 
\begin{equation}\label{eq_bound1_ob}
	\int_\Omega \abs{(\partial_t u) (X\oline{u})} \de V \le \frac{L}{\sqrt2} E(0)
\end{equation}
holds here as well. We can also get the same bound on the first term -- the argument just requires a few sign changes. In this case, the function $F$ in \eqref{eq_F} must now be defined as 
\begin{equation}
	F(x) = \int_{\frac{a_1}{\ell} x}^{\frac{a_2}{\ell} x} \abs{u(t_0,x,y)}^2 \de y,
\end{equation}
and we compute 
\begin{align}
\begin{split}
	F'(x) &= \frac{a_2}{\ell} \Abs{u\paren{t_0,x,\frac{a_2}{\ell} x}}^2 - \frac{a_1}{\ell} \Abs{u\paren{t_0,x,\frac{a_1}{\ell} x}}^2 + 2 \re \int_{\frac{a_1}{\ell} x}^{\frac{a_2}{\ell} x} u \partial_x u \de y \\
		&\le \gamma \int_{\frac{a_1}{\ell} x}^{\frac{a_2}{\ell} x} \abs{u}^2 \de y 
			+ \frac 1\gamma \int_{\frac{a_1}{\ell} x}^{\frac{a_2}{\ell} x} \abs{\partial_x u}^2 \de y \\
		&= \gamma F(x) + \frac 1\gamma G(x)
\end{split}
\end{align}
where the function $G$ is now 
\begin{equation}
G(x) = \int_{\frac{a_1}{\ell} x}^{\frac{a_2}{\ell} x} \abs{\partial_x u}^2 \de y. 
\end{equation}
From this point the same sequence of steps applies, and we obtain 
\begin{equation}
	\Abs{ \int_\Omega \partial_t u \oline{u} \de V } \le \frac12 L \sqrt{e-1} E(0)
\end{equation}
as in \eqref{eq_bound2}. This and \eqref{eq_bound1_ob} imply 
\begin{equation}\label{ineq_with_rhs_ob}
	\Abs{\int_0^T\int_{\partial\Omega} (Xu)(\partial_\nu\overline{u}) \,dS \,dt} 
		= TE(0)\paren{1 + \mc{O}\paren{\frac{L}{T}}}
\end{equation}
mirroring \eqref{ineq_with_rhs}.

Finally, we obtain the Neumann data from the term on the left hand side of \eqref{eq_comm_equality} via similar steps. On $A$, we again have $\partial_y u = 0$, $x=\ell$, and $\partial_\nu u = \partial_x u$, so $Xu = \ell \partial_\nu u$ and 
\begin{align}\label{eq_A_data_ob}
\begin{split}
	\int_0^T\int_A (Xu)(\partial_\nu\overline{u}) \,dS \,dt 
		&= \int_0^T\int_A \ell (\partial_\nu u)(\partial_\nu\overline{u}) \,dS \,dt \\
		&= \int_0^T\int_A \ell \abs{\partial_\nu u}^2 \,dS \,dt.
\end{split}
\end{align}

On side $B$, we have $y = \frac{a_1}{\ell} x$, and the unit tangent vector is 
$\tau = \paren{\frac{\ell}{b}, \frac{a_1}{b}}$. Since the tangential derivative vanishes, we have 
\begin{equation}\label{eq_B_relation_ob}
	\tau \cdot \nabla u = \frac{\ell}{b} \partial_x u + \frac{a_1}{b} \partial_y u = 0.
\end{equation}
The unit normal vector on $B$ is $\nu = \paren{\frac{a_1}{b}, -\frac{\ell}{b}}$. Using this and \eqref{eq_B_relation_ob}, we obtain 
\begin{align}
\begin{split}
	(Xu)(\partial_\nu\oline{u}) &= \paren{x\partial_x u + \frac{a_1}{\ell}x\partial_y u} \paren{\frac{a_1}{b}\partial_x\oline{u} - \frac{\ell}{b}\partial_y\oline{u}} \\
		&= x\paren{\frac{a_1}{b}\partial_xu\partial_x\oline{u} + \frac{a_1^2}{b\ell} \partial_yu\partial_x\oline{u} 
			- \frac{\ell}{b}\partial_xu\partial_y\oline{u} - \frac{a_1}{b}\partial_yu\partial_y\oline{u}} \\
		&= x\frac{a_1}{\ell}\partial_x\oline{u}\paren{\frac{\ell}{b} \partial_x u + \frac{a_1}{b} \partial_y u}
			- x\partial_y\oline{u}\paren{\frac{\ell}{b} \partial_x u + \frac{a_1}{b} \partial_y u} \\
		&= 0,
\end{split}
\end{align}
and consequently 
\begin{equation}\label{eq_B_data_ob}
	\int_0^T\int_B (Xu)(\partial_\nu\overline{u}) \,dS \,dt = 0.
\end{equation}

In a similar fashion, on side $C$ we have $y = \frac{a_1+a_2}{c}x$, and the unit tangent vector is 
$\tau = \paren{\frac{\ell}{c}, \frac{a_1+a_2}{c}}$. Since the tangential derivative vanishes, we have 
\begin{equation}\label{eq_C_relation_ob}
	\tau \cdot \nabla u = \frac{\ell}{c} \partial_x u + \frac{a_1+a_2}{c} \partial_y u = 0.
\end{equation}
The unit normal vector on $C$ is $\nu = \paren{-\frac{a_1+a_2}{c}, \frac{\ell}{c}}$. Using this and \eqref{eq_C_relation_ob}, 
we obtain 
\begin{align}
\begin{split}
	(Xu)&(\partial_\nu\oline{u}) \\ 
		&= \paren{x\partial_x u + \frac{a_1+a_2}{\ell}x\partial_y u} 
			\paren{-\frac{a_1+a_2}{c}\partial_x\oline{u} + \frac{\ell}{c}\partial_y\oline{u}} \\
		&= x\paren{-\frac{a_1+a_2}{c}\partial_xu\partial_x\oline{u} 
			- \frac{(a_1+a_2)^2}{c\ell} \partial_yu\partial_x\oline{u} 
			+ \frac{\ell}{c}\partial_xu\partial_y\oline{u} 
			+ \frac{a_1+a_2}{c}\partial_yu\partial_y\oline{u}} \\
		&= -x\frac{a_1+a_2}{\ell}\partial_x\oline{u}\paren{\frac{\ell}{c} \partial_x u 
			+ \frac{a_1+a_2}{c} \partial_y u}
			+ x\partial_y\oline{u}\paren{\frac{\ell}{c} \partial_x u + \frac{a_1+a_2}{c} \partial_y u} \\
		&= 0,
\end{split}
\end{align}
and consequently 
\begin{equation}\label{eq_C_data_ob}
	\int_0^T\int_C (Xu)(\partial_\nu\overline{u}) \,dS \,dt = 0.
\end{equation}

Thus, from \eqref{eq_B_data_ob} and \eqref{eq_C_data_ob} we have 
\begin{equation}
	\int_0^T\int_{\partial\Omega} (Xu)(\partial_\nu\overline{u}) \,dS \,dt = \int_0^T\int_A \ell \abs{\partial_\nu u}^2 \,dS \,dt
\end{equation}
and \eqref{ineq_with_rhs_ob} gives  
\begin{equation}
	\int_0^T\int_A \ell \abs{\partial_\nu u}^2 \,dS \,dt = TE(0) \paren{1 + \mc{O}\paren{\frac{L}{T}}},
\end{equation}
implying 
\begin{equation}
	\int_0^T\int_A \abs{\partial_\nu u}^2 \,dS \,dt = \frac{T}{\ell}E(0) \paren{1 + \mc{O}\paren{\frac{L}{T}}} 
\end{equation}
which completes the proof of Theorem \ref{thm_main}.

\section{Failure of the Result on Square Domains}\label{ex_square}

We again draw inspiration from \cite{Chr2DTriangles} to construct an example to demonstrate that our main 
result fails for general polygons -- in particular, for squares. Let $\Omega = [0, 2\pi]^2$, and consider 
\begin{equation}
	\phi(x,y) = \frac1\pi \sin(x)\sin(ny)
\end{equation}
for some integer $n>0$. Note that $\phi$ is a Dirichlet eigenfunction satisfying 
\begin{equation}
	\begin{cases}
		(-\Delta - (1+n^2)) \phi = 0 \\
		\phi |_{\partial\Omega} = 0,
	\end{cases}
\end{equation}
which is normalized: $\norm[L^2(\Omega)]{\phi} = 1$. Define now 
\begin{equation}\label{square_efun}
	u(t,x,y) = \sin(t\sqrt{1+n^2})\phi(x,y),
\end{equation}
which satisfies 
\begin{equation}
	\begin{cases}
		(\partial_t^2-\Delta)u = 0 \\
		u |_{\partial\Omega} = 0 \\
		u(0,x,y) = 0,\ u_t(0,x,y) = \sqrt{1+n^2}\phi(x,y).
	\end{cases}
\end{equation}
The energy is then 
\begin{align}
\begin{split}
	E(0) &= \int_\Omega \abs{\partial_t u}^2 + \abs{\nabla u}^2 \de V \\
		&= \int_\Omega (1+n^2)\abs{\phi}^2 + 0 \de V \\
		&= 1+n^2.
\end{split}
\end{align}
Along the right edge $\{2\pi\} \by [0,2\pi]$, where $\partial_\nu u = \partial_x u$, we then have 
\begin{align}
\begin{split}
	\int_0^T\int_0^{2\pi} & \abs{\partial_\nu u}^2
        \Big\vert_{x=2\pi} \de y \de t \\ 
		&= \int_0^T\int_0^{2\pi} \frac{1}{\pi^2} \abs{ \sin(\sqrt{1+n^2} t) \cos(2\pi) \sin(ny) }^2 \de y \de t \\
		&= \frac{1}{\pi^2} \int_0^T\int_0^{2\pi} \sin^2(\sqrt{1+n^2} t) \sin^2(ny) \de y \de t \\
		&= \frac{1}{4\pi^2} \int_0^T\int_0^{2\pi} (1-\cos(2\sqrt{1+n^2} t)) (1-\cos(2ny)) \de y \de t \\
		&= \frac{1}{4\pi^2} \paren{ 2\pi T - \frac{2\pi}{2\sqrt{1+n^2}} \sin(4\pi T\sqrt{1+n^2}) } \\
		&= \frac{T}{2\pi} - \frac{1}{4\pi\sqrt{1+n^2}}
        \sin(4\pi T\sqrt{1+n^2}) \\
        & = \left(\frac{T}{2\pi(1 + n^2)} - \frac{1}{4\pi(1+n^2)^{3/2}}
        \sin(4\pi T\sqrt{1+n^2}) \right)E(0).
\end{split}
\end{align}
From this, we see that the analog of Theorem \ref{thm_main} fails to hold in this context, as an observability inequality cannot be established: there is no $T>0$ and a constant 
$C_T$ depending only on T and $\Omega$ such that 
\begin{equation}
	\int_0^T\int_0^{2\pi} \abs{\partial_\nu u}^2 \Big\vert_{x=2\pi} \de y \de t \ge C_T E(0)
\end{equation}
for all solutions $u$, as for any $C_T$, a suitably large $n$ may be
chosen so that \eqref{square_efun}  serves as a counterexample.

%%%%%%%%%%%%%%%%%%%%%%%%%%%%%%%%%%%%%%%%%%%%%%%%%%%%%%%%%%%%%
%% APPENDICES
%%%%%%%%%%%%%%%%%%%%%%%%%%%%%%%%%%%%%%%%%%%%%%%%%%%%%%%%%%%%%

% Bibliography
\bibliographystyle{alpha}
\bibliography{triangle-wave-bib}

% Add bibliography to contents page
%\addcontentsline{toc}{chapter}{Bibliography} %'Bibliography' into toc

\end{document}